\documentclass[11pt,a4paper,reqno]{amsart}
\usepackage[top=1.15in, bottom=1.15in, left=1.17in, right=1.17in]{geometry}
\usepackage[colorlinks=true, pdfstartview=FitV, linkcolor=blue, citecolor=blue, urlcolor=blue, breaklinks=true]{hyperref}

\usepackage{hyphenat}
\usepackage{amssymb}
\usepackage{mathrsfs}
\usepackage{todonotes}
\usepackage{mathtools}
\usepackage{amsthm}
\usepackage{amscd}
\usepackage{todonotes}
\usepackage{enumitem}
\usepackage{comment}
\usepackage{graphicx}
\usepackage{dsfont}
\usepackage{nameref, cleveref}
\usepackage{fancyhdr, color,tikz-cd}
\usepackage[table]{xcolor}
\usepackage{subcaption}

\makeatletter
\newcommand*{\@old@slash}{}\let\@old@slash\slash
\def\slash{\relax\ifmmode\delimiter"502F30E\mathopen{}\else\@old@slash\fi}
\makeatother

\newcommand{\g}{\mathfrak{g}}
\newcommand{\gl}{\mathfrak{gl}}

\newcommand{\rep}{\mathrm{Rep}}
\newcommand{\e}{\varepsilon}
\newcommand{\proj}{\pi_{(\underline a;\underline c),(\underline b;\underline d)}}

\newcommand{\C}{\mathbb{C}}

\newcommand{\F}{\mathcal{F}}
\newcommand{\Z}{\mathbb{Z}}
\newcommand{\N}{\mathbb{N}}

\renewcommand{\o}{\circ}

\renewcommand{\a}{\alpha}

\renewcommand{\d}{\delta}
\renewcommand{\b}{\mathfrak{b}}

\newcommand{\h}{\mathfrak{h}}

\newcommand{\m}{\mathbf{m}}

\renewcommand{\sl}{\mathfrak{sl}}

\numberwithin{equation}{section}

\renewcommand{\hom}{\mathsf{Hom}}

\renewcommand{\a}{\underline a}
\renewcommand{\b}{\underline b}
\newcommand{\W}{\mathcal{W}}
\newcommand{\RW}{\mathcal{RW}}

\newtheorem{thm}{Theorem}[section]
\newtheorem*{thm*}{Theorem}
\newtheorem{corollary}[thm]{Corollary}
\newtheorem{lemma}[thm]{Lemma}
\newtheorem{proposition}[thm]{Proposition}

\newtheorem{algorithm}[thm]{Algorithm}

\theoremstyle{definition}
\newtheorem{definition}[thm]{Definition}

\theoremstyle{remark}
\newtheorem{remark}[thm]{Remark}

\crefname{lemma}{lemma}{lemmas}
\Crefname{Lemma}{Lemma}{Lemmas}

\crefname{thm}{theorem}{theorems}
\Crefname{thm}{Theorem}{Theorems}

\title{Web bases for the Lie superalgebra $\mathfrak{gl}(1|1)$}

\author[Enugandla]{Pranav Enugandla}
\address[Enugandla]{University of California, Berkeley, CA, United States of America}
\email{shreepranav\_varma@berkeley.edu}

\date{\today}

\begin{document}

\begin{abstract}
Queffelec--Sartori \cite{queffelec-sartori} and Tubbenhauer--Vaz--Wedrich \cite{tubbenhauer-vaz-wedrich} give a diagrammatic web category via generators and relations to describe the subcategory of $\rep(\gl(m|n))$ generated by tensor products of exterior powers of the defining representation. We construct the first explicit bases for homomorphism spaces in the web category for $\gl(1|1)$, giving an alternative diagrammatic presentation of this web category, and explore connections to the combinatorics of planar networks. Our basis comes equipped with an algorithm to decompose any web into a linear combination of basis webs, giving an effective way to compute in the $\gl(1|1)$-skein theory.
\end{abstract}

\maketitle
\section{Introduction}

\subsection{Webs for classical Lie algebras}
\emph{Webs} were introduced by Kuperberg in \cite{Kuperberg} as a diagrammatic calculus for effectively computing quantum link invariants, via a diagrammatic presentation for $U_q(\sl_3)$ invariant homomorphism between tensor products of fundamental representations. More precisely, if $V_q^{\a} = \bigotimes_{i=1}^k V_q(\omega_{a_i})$ and $V_q^{\b} = \bigotimes_{j=1}^\ell V_q(\omega_{b_j})$, then \cite{Kuperberg} describes a diagrammatic basis for the space of $U_q(\sl_3)$ invariant homomorphisms:
\[\hom_{U_q(\sl_3)}(V_q^{\a}, V_q^{\b}).\]

This basis is described by \emph{webs}, which are planar trivalent bipartite graphs drawn between two horizontal lines: the bottom boundary vertices labeled $a_1,\dots, a_k$ and the top boundary vertices labeled $b_1,\dots, b_\ell$ from left to right. The morphisms denoted by webs usually have many relations between each other, and a key problem is to find a \emph{web basis}, i.e., a subset of the web invariants forming a basis for the morphism space. In his seminal paper on webs, Kuperberg finds the \emph{non-elliptic web basis} in the case of $U_q(\sl_3)$, as well as the other rank $2$ Lie algebras, by forbidding certain local configurations. The natural problem then was to extend this theory to higher rank $\sl_n$.

The major breakthrough in this direction was by Cautis, Kamnitzer, and Morrison in \cite{Cautis-Kamnitzer-Morrison}, where they give a diagrammatic presentation for these morphism spaces for arbitrary rank $\sl_n$, in particular describing a generating set for all the relations between web invariants. Since then, there have been many approaches to finding web bases in higher rank, including \cite{Fontaine}, \cite{Elias}, and the first rotation invariant web basis for $U_q(\sl_4)$ in \cite{GPPSS-sl4}.

Since their invention, webs have found connections and applications to quantum link invariants \cite{Khovanov}, skein modules \cite{Le.Sikora}, dimer models \cite{Douglas-Kenyon-Shi}, dynamical algebraic combinatorics \cite{Petersen-Pylyavskyy-Rhoades}, and cluster algebras \cite{Fomin-Pylyavskyy-advances}. 

\subsection{Webs in the context of tangle invariants} One of the earliest applications of webs was to provide an efficient algorithm for computing link invariants and tangle invariants. Reshetikhin and Turaev \cite{Reshetikhin-Turaev} define an invariant of oriented framed tangles for each ribbon Hopf algebra, usually taken to be the quantization $U_q(\g)$ of the universal enveloping algebra of a Lie algebra $\g$. When $\g=\mathfrak{sl}_2$, this recovers the familiar Jones polynomial of a link.

Kuperberg's webs give an efficient algorithm for computing these tangle invariants for $\g$ of rank $2$ by resolving a tangle projection into the non-elliptic web basis (see \cite[Section 8.3]{Kuperberg}). More generally, one can consider stated $\g$-skein algebras of surfaces, and web bases give an effective way to compute in these skein theories. 

Another classical invariant of links is the Alexander polynomial. However, it is not very straightforward to interpret it within the Reshetikhin--Turaev framework. One of the main difficulties, as noted in \cite{grant16}, is that the Reshetikhin--Turaev invariant associates to the unknot the quantum dimension of the representation coloring it, whereas the Alexander polynomial is $0$ for a split link. This urges us to look for representation theories that admit objects of dimension $0$ in the appropriate categorical sense. The representation theory of the quantized superalgebra $U_q(\gl(1|1))$ is perhaps the simplest such theory. This perspective is thoroughly investigated in \cite{Sartori2013TheAP}, wherein it is shown how to recover the Alexander polynomial from the Reshetikhin--Turaev tangle invariant coming from $U_q(\gl(1|1))$. A natural extension is to study the $U_q(\gl(1|1))$ skein theory on more general surfaces, which would benefit from a good web basis.

\subsection{Webs for \texorpdfstring{$\gl(m|n)$}{gl(m|n)}}

The above motivates us to extend web bases to the general linear Lie superalgebra $\gl(m|n)$ and its quantum version $U_q(\gl(m|n))$. These are $\Z/2\Z$ graded analogues of $\gl_n$ and $U_q(\gl_n)$. For an introduction to the representation theory of Lie superalgebras, see \cite{cheng-wang} or \cite{musson}. Web diagrammatics for the special case of $U_q(\gl(1|1))$ were investigated in \cite{grant16}, \cite{Sartori2013CategorificationOT}. This work was built upon in \cite{queffelec-sartori}, \cite{queffelec2015noteglmnlinkinvariants}, and later in \cite{tubbenhauer-vaz-wedrich}, wherein a diagrammatic presentation of a certain subcategory of $\mathrm{Rep}(U_q(\gl(m|n)))$ is given in terms of generators and relations. These articles also explore the connection between the Reshetikhin--Turaev invariants coming from $\rep(U_q(\gl(m|n)))$ and the HOMFLY-PT polynomial, showing in particular that these invariants only depend on the integer $m-n$. 

Henceforth, we let $\g = \gl(1|1)$, and explain the non-quantum story for simplicity. The representation category in question is $\mathrm{Fund}(\g)\subset \rep(\g)$, the category monoidally generated by the exterior powers $\bigwedge^nV$ of the defining representation $V\simeq \C^{1|1}$ of $\g$. Even though the category $\rep(\g)$ of finite-dimensional representations is \emph{not} semisimple, $\mathrm{Fund}(\g)$ lives inside the semisimple category $\rep^+(\g)\subset \rep(\g)$ of \emph{polynomial representations}. A major difference from the classical $\gl_n$ picture is that arbitrarily high exterior powers of representations are non-vanishing in the super case. In terms of webs, this means that edges can have arbitrarily high multiplicities.

We start with the \emph{free web category} $\infty\mathsf{-Web}^f$ of trivalent planar graphs, all of whose edges are oriented upwards and have multiplicites such that the algebraic sum at each internal vertex is $0$. We draw all our webs between two horizontal boundaries, to be read as morphisms from the bottom to the top. There is a realization functor 
\[\alpha: \infty\mathsf{-Web}^f\longrightarrow \mathrm{Fund}(\g).\]
It is a standard fact that $\alpha$ is full. In \cite{tubbenhauer-vaz-wedrich}, an explicit description of the kernel of $\alpha$ is provided by imposing relations on $\infty\mathsf{-Web}^f$. However, the problem of finding good \emph{web bases} for the morphism spaces still remains.

\subsection{Main results}

Our key insight is that under the realization functor $\alpha$, the webs in \Cref{fig:reduction} can be resolved (see (\ref{eq:square_rel1}) and (\ref{eq:zipper_rel})):

\begin{figure}[ht]
    \centering
    \includegraphics[height=3cm]{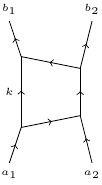}
    \hspace{3cm}
    \includegraphics[height=3cm]{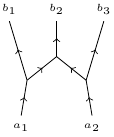}
    \caption{Webs that can be resolved into simpler webs. We call these \emph{squares} (on the left) and \emph{zippers} (on the right).}
    \label{fig:reduction}
\end{figure}

This motivates the definition of \emph{reduced webs} to be those that do not contain any bigons, squares, or zippers. The resolving relations show that reduced webs span the homomorphism spaces in question. The main result of the article is that they in fact form a basis.

\begin{thm}
    Reduced webs give bases for morphism spaces in $\mathrm{Fund}(\g)$.
\end{thm}

The relations resolving squares and zippers give an algorithm to decompose a web into basis elements: 

\begin{algorithm}[Reduction algorithm] Let $W$ be any web. The following steps reduce $W$ to a linear combination of reduced webs.
\begin{enumerate}
    \item Resolve all bigons using (\ref{eq:bigon}).
    \item Fix a face of $W$. 
    \begin{itemize}
        \item If the chosen face has $\ge6$ edges, then resolve using (\ref{eq:zipper_rel}), to reduce the number of edges in the chosen face.
        \item If the chosen face is a square, resolve it using (\ref{eq:square_rel1}).
    \end{itemize}
    Repeat for all other faces of $W$.
    \item Resolve all zippers using (\ref{eq:zipper_rel}).
\end{enumerate}
\end{algorithm}

This gives an effective way of computing in the $\gl(1|1)$-skein theory. As a corollary, we obtain an alternative set of relations on $\infty\mathsf{-Web}^f$, such that $\alpha$ descends to an equivalence of categories between the quotient and $\mathrm{Fund}(\g)$.

Moreover, we extend this theory to \emph{clasped webs}, to access the larger category $\rep^+(\g)$ of all polynomial representations of $\gl(1|1)$. The idea is to embed any irreducible representation in $\rep^+(\g)$ as the Cartan component of a representation of $\mathrm{Fund}(\g)$. More generally, given a tensor product of irreducibles $L(c_1|p_1)\otimes \cdots\otimes L(c_k|p_k)$ in $\rep^+(\g)$, there is a natural way to embed it inside some $V^{\a}\in \mathrm{Fund}(\g)$. Picking a similar embedding $L(d_1|q_1)\otimes \cdots \otimes L(d_\ell|q_\ell)\subset V^{\b}$, we see that
\[\hom_\g\big(L(c_1|p_1)\otimes \cdots\otimes L(c_k|p_k), L(d_1|q_1)\otimes \cdots \otimes L(d_\ell|q_\ell)\big) \subset \hom_\g(V^{\a}, V^{\b}).\]
We give a web basis for this subspace using the notion of \emph{non-convexity} for $\gl(1|1)$ webs, motivated by \cite{Kuperberg} in the $\sl_3$ case and the local description in \cite{enugandla-gaetz} in the $\sl_4$ case.
\begin{thm}
    In the above situation, the non-convex webs in the web basis for $\hom_\g(V^{\a}, V^{\b})$ give a basis for $\hom_\g\big(L(c_1|p_1)\otimes \cdots\otimes L(c_k|p_k), L(d_1|q_1)\otimes \cdots \otimes L(d_\ell|q_\ell)\big)$. 
\end{thm}

\subsection{Outline}

In \Cref{sec:prelim} we recall background on Lie superalgebras following \cite{cheng-wang}, and the diagrammatic web category from \cite{tubbenhauer-vaz-wedrich}. In \Cref{sec:flows} we provide a combinatorial description of the monomial expansion of web invariants using flows. In \Cref{sec:rels} we precisely state and prove the main result about the basis of reduced webs. In \Cref{sec:clasped-webs} we extend our web basis to the clasped setting, where tensor products of more general irreducible $\gl(1|1)$ representations are considered.

\section{Preliminaries}
\label{sec:prelim}

\subsection{Notation}

Throughout the article, we work over the ground field $\C$. We let $\Z_{\ge0}$ denote the set of nonnegative integers, and $\N = \Z_{>0}$ the set of positive integers. For $n\in \N, k\in \Z_{\ge0}$, we denote by $[n]$ the set $\{1,\dots, n\}$, and by $\binom{[n]}{k}$ the set of $k$ element subsets of $[n]$. For a tuple $\a = (a_1,\dots,a_k)\in \N^k$, we let $|\a|$ denote the sum $\sum_{i=1}^k a_i$.

For a $\Z/2\Z$ graded vector space $V = V_{\bar0}\oplus V_{\bar 1}$ and a homogeneous vector $v\in V$, we denote by $p(v)\in\Z/2\Z$ the \emph{parity} of $v$, defined as
\[p(v):=\begin{cases} 
\bar0 &\text{if }v\in V_{\bar0},\\
\bar1 &\text{if }v\in V_{\bar1}.
\end{cases}\]
We call a homogeneous element $v\in V$ \emph{even} (resp. \emph{odd}) if $p(v)=\bar0$ (resp. $p(v)=\bar1$).

\subsection{Lie superalgebra \texorpdfstring{$\gl(1|1)$}{gl(1|1)}} 

We use \cite{cheng-wang} as our main reference for the representation theory of the Lie superalgebra $\g :=\gl(1|1)$. 

Let $V$ denote the $1|1$-dimensional superspace $\C^{1|1}$, with one even dimension spanned by $x$ and one odd dimension spanned by $y$, i.e.,
\[V_{\bar0} = \C x = \C\begin{pmatrix}
1\\ 0\end{pmatrix},\quad V_{\bar1} = \C y=\C\begin{pmatrix}
0\\ 1\end{pmatrix}.\]
The Lie superalgebra of endomorphisms of $V$ is $\g:=\gl(V)\simeq \gl(1|1)$. Explicitly, $\gl(1|1)$ is the space of $2\times 2$ matrices with even part
\[\gl(1|1)_{\bar0} = \mathrm{span}_\C\left\{
h_1 =\left(
\begin{array}{c|c}
    1 &0  \\
    \hline
    0 & 0
\end{array}
\right),\ 
h_2 =\left(
\begin{array}{c|c}
    0 & 0  \\
    \hline
    0 & 1
\end{array}
\right)
\right\},\]
and odd part
\[\gl(1|1)_{\bar1} = \mathrm{span}_\C\left\{
e =\left(
\begin{array}{c|c}
    0 & 1  \\
    \hline
    0 & 0
\end{array}
\right),\ 
f =\left(
\begin{array}{c|c}
    0 & 0  \\
    \hline
    1 & 0
\end{array}
\right)
\right\}.\]

On homogeneous elements $A$ and $B$, the bracket of the Lie superalgebra structure is given by
\[[A, B] := AB - (-1)^{p(A)p(B)} BA.\]

With respect to the basis chosen, the Lie bracket can be described as follows:
\[[h_1, e] = e,\quad [h_1, f] = -f,\]
\[[h_2, e] = -e,\quad [h_2, f] = f,\]
\[[h_1, h_2] = [e,e] = [f,f] = 0,\]
\[[e, f] = h_1+h_2 =: h.\]

Note that $h\in Z(\g)$. We define the \emph{Cartan subalgebra} $\h$ to be the subalgebra of diagonal matrices, i.e.,
\[\h := \mathrm{span}_\C\{h_1, h_2\}\]
We define $\{\e,\d\}\subset \h^*$ to be the basis dual to $\{h_1, h_2\}\subset \h$. Then the \emph{weight lattice} is defined to be
\[\Lambda:= \Z\e\oplus \Z \d.\]
For ease of notation, we denote the weight $a\e+b\d\in \Lambda$ by $(a|b)$.
Our choice of Borel subalgebra $\mathfrak b$ will be
\[\mathfrak b := \h \oplus \C e.\]
This corresponds to the positive root being $\e - \d = (1|-1)\in \Lambda$. Given any $\lambda\in \Lambda$, we get a highest weight irreducible representation $L(\lambda)$, defined as the quotient of the Verma module $M(\lambda):=\mathcal {U}\g\otimes_{\mathcal{U}\mathfrak{b}} \C_{\lambda}$ by its maximal submodule.

\subsection{Polynomial representations of \texorpdfstring{$\gl(1|1)$}{gl(1|1)}}

The \emph{defining representation} of $\gl(1|1)$ is $V\simeq \C^{1|1}\simeq L(1|0)$, where the action is via matrix multiplication. Given two $\gl(1|1)$ representations $W_1$ and $W_2$, we define the structure of a representation of $W_1\otimes W_2$ via
\[x\cdot (w_1\otimes w_2) = x\cdot w_1\otimes w_2 + (-1)^{p(x)p(w_1)}w_1\otimes x\cdot w_2,\]
where $x\in \gl(1|1), w_1\in W_1, w_2\in W_2$ are homogeneous.

This gives $T^dV:= V^{\otimes d}$ the structure of a $\gl(1|1)$ representation. It is known (for example, by Schur--Sergeev duality \cite[Section 3.2]{cheng-wang}) that $T^dV$ is semisimple.

\begin{definition}
    The category $\rep^+(\g)$ of \emph{polynomial representations} is the category of representations of $\g$ which occur as direct sums of direct summands of $T^dV$.
\end{definition}

\begin{remark}
    Polynomial representations are referred to as \emph{covariant representations} in some literature.
\end{remark}

$\rep^+(\g)$ is a semisimple abelian category that admits a symmetric monoidal $\otimes$-structure, i.e., it is a semisimple symmetric tensor category. The irreducibles in $\rep^+(\g)$ are all of the form $L(\lambda)$ for $\lambda\in \Lambda^{\mathrm{pol}}$, where
\[\Lambda^{\mathrm{pol}} = \{(a|b): a, b\in \Z_{\ge 0}\text{ such that } a=0 \implies b=0\}\]

Of particular interest to us will be the exterior powers $\bigwedge^dV\in \rep^+(\g)$. Formally, $\bigwedge^\bullet V$ is the quotient of $T^\bullet V$ by the two sided ideal generated by $v\otimes w + (-1)^{p(v)p(w)}w\otimes v$ for homogeneous $v,w\in V$. $\bigwedge^d V$ is then the degree $d$ component of $\bigwedge^\bullet V$.

An easy computation shows that the basis $\{x,y\}$ for $V$ induces the basis $\{x\wedge y^{d-1}, y^d\}$ for $\bigwedge^d V$. Note that $x\wedge y^{d-1}$ has weight $(1|d-1)$, and $y^d$ has weight $(0|d)$. One can check that $\bigwedge^d V$ is irreducible; more precisely:
\[\bigwedge^d V \simeq L(1|d-1).\]

\begin{definition}
    $\mathrm{Fund}(\g)$ is the full subcategory of $\rep^+(\g)$ whose objects are of the form 
    \[V^{\a} := \bigwedge^{a_1}V\otimes \cdots \otimes \bigwedge^{a_k}V,\quad \a = (a_1, \dots, a_k)\in \N^k.\]
\end{definition}

We stress that $\mathrm{Fund}(\g)$ is \emph{not} an abelian category (not even an additive category), and is merely a $\C$-linear symmetric monoidal category.

\subsection{The diagrammatic web category}

A diagrammatic presentation for $\mathrm{Fund}(\g)$ is described in \cite{tubbenhauer-vaz-wedrich} using a super version of quantum skew Howe duality. This generalizes earlier presentations by \cite{grant16}, \cite{Sartori2013CategorificationOT}, and \cite{queffelec-sartori}.

\begin{definition}[see \cite{tubbenhauer-vaz-wedrich}]
     $\infty\mathsf{-Web}^f$ is the $\C$-linear monoidal category of planar directed trivalent graphs (up to relative height preserving isotopy) with all edges oriented upwards, such that the algebraic sum of edge multiplicities at each internal vertex is $0$. In other words, $\infty\mathsf{-Web}^f$ is described by

    \begin{itemize}
        \item \emph{Objects:} Tuples of positive integers,
        \item \emph{Morphisms:} Monoidally generated by \emph{simple merge} and \emph{simple split} vertices, as pictured in \Cref{fig:merge_split}.
        \begin{figure}
            \centering
            \includegraphics[height=2cm]{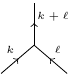}
            \hspace{3cm}
            \includegraphics[height=2cm]{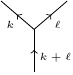}
            \caption{A simple merge (left) and a simple split (right) vertex.}
            \label{fig:merge_split}
        \end{figure}
    \end{itemize}
\end{definition}

\begin{definition}
    $\infty\mathsf{-Web}$ is the quotient of $\infty\mathsf{-Web}^f$ by the following relations (including the horizontal reflection of the last relation):
    
    \begin{equation}\label{eq:assoc}
    \vcenter{\hbox{\includegraphics[width=0.5\linewidth]{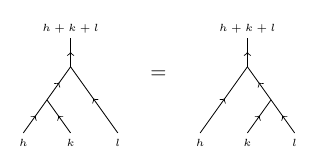}}}
    \end{equation}
    \begin{equation}\label{eq:coassoc}
    \vcenter{\hbox{\includegraphics[width=0.5\linewidth]{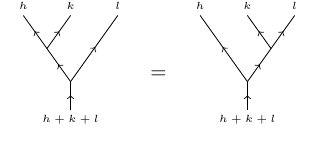}}}
    \end{equation}
    \begin{equation}\label{eq:bigon}
    \vcenter{\hbox{\includegraphics[width=0.5\linewidth]{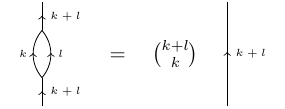}}}
    \end{equation}
    \begin{equation}
    \vcenter{\hbox{\includegraphics[width=0.9\linewidth]{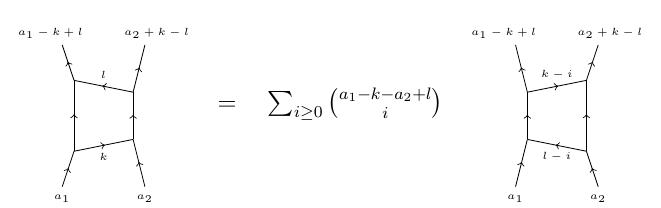}}}        
    \end{equation}
\end{definition}

\begin{remark}\label{rem:simplify}
    The relation (\ref{eq:assoc}) (resp. (\ref{eq:coassoc})) allows us to use multivalent merges (resp. multivalent splits) to denote the composition of trivalent merges (resp. trivalent splits). Using multivalent vertices, we can, and henceforth will, assume that the morphisms in $\infty\mathsf{-Web}$ are \emph{merge-split bipartite}, i.e., no two merges or two splits are adjacent. Additionally, we omit the orientations of edges, since they all point upwards.
\end{remark}

There is a realization functor $\alpha: \infty\mathsf{-Web}\to \mathrm{Fund}(\gl(m|n))$ for any $\gl(m|n)$, defined on objects as $\a=(a_1,\dots, a_k)$ mapping to $V^{\a}$, and on morphisms by:

\begin{align*}
\vcenter{\hbox{\includegraphics[height=2cm]{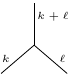}}}\quad&:\quad \bigwedge^kV\otimes \bigwedge^\ell V\longrightarrow \bigwedge^{k+\ell}V\\
& \hspace{1.4cm} v\otimes w\quad\longmapsto v\wedge w
\end{align*}

\begin{align*}
\vcenter{\hbox{\includegraphics[height=2cm]{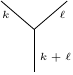}}}\quad&:\quad \bigwedge^{k+\ell}V \longrightarrow \bigwedge^k V\otimes \bigwedge^\ell V\\
& \hspace{2mm} v_1\wedge\cdots\wedge v_{k+\ell}\quad\longmapsto \sum_{S\in \binom{[k+\ell]}{k}}\mathrm{sgn}(S)\ v_S\otimes v_{[k+\ell]\setminus S}
\end{align*}

where if $S = \{i_1<\cdots<i_k\}$ and $[k+\ell]\setminus S = \{j_1<\cdots<j_\ell\}$, then $v_S:= v_{i_1}\wedge\cdots\wedge v_{i_k}$, $v_{[k+\ell]\setminus S}:= v_{j_1}\wedge\cdots\wedge v_{j_\ell}$, and $\mathrm{sgn}(S):= (-1)^{(i_1-1)+\cdots +(i_k-k)}$.

It is known that $\alpha$ is essentially surjective and full. The results of \cite[Section 5]{tubbenhauer-vaz-wedrich} give a description of the kernel of $\alpha$. We henceforth specialize to $\g = \gl(m|n) =\gl(1|1)$. 

\begin{definition}
    $1|1\mathsf{-Web}$ is the quotient of $\infty\mathsf{-Web}$ by the relation (and its horizontal reflection) below:
    \begin{equation}\label{eq:funny}
        a_1a_2(a_1-1)(a_2+1)
        \vcenter{\hbox{\includegraphics[height=2.7cm]{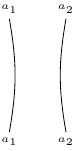}}}
        -\quad 2(a_1-1)(a_2+1)
        \vcenter{\hbox{\includegraphics[height=2.7cm]{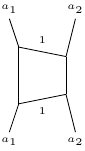}}}
        + \quad 4
        \vcenter{\hbox{\includegraphics[height=2.7cm]{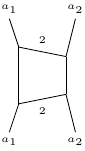}}}
        = 0.
    \end{equation}
\end{definition}

The content of \cite[Theorem 5.11]{tubbenhauer-vaz-wedrich} along with \cite[Example 5.3]{tubbenhauer-vaz-wedrich} shows the following:

\begin{thm}[\cite{tubbenhauer-vaz-wedrich}]\label{thm:tvw}
    $\alpha$ descends to an equivalence of categories
    \[\alpha: {1|1}\mathsf{-Web}\to \mathrm{Fund}(\g)\]
\end{thm}

The main goal of this article is to find an alternative set of relations in $\infty\mathsf{-Web}^f$ which give a diagrammatic presentation of $\mathrm{Fund}(\g)$. These relations will have an obvious directionality, in the sense that these relations give a way of resolving forbidden configurations. This allows us to find bases for morphism spaces, instead of just a generators and relations presentation.

\section{Flows and \texorpdfstring{$\gl(1|1)$}{gl(1|1)} webs}
\label{sec:flows}

Our goal is to find an alternative diagrammatic presentation of $\mathrm{Fund}(\g)$, with the property that the relations have an obvious directionality. This will allow us to find diagrammatic bases for the morphism spaces. 

\subsection{Dimensions of homomorphism spaces}
Throughout this section $\g = \mathfrak{gl}(1|1)$. Fix tuples of nonnegative integers $\underline a = (a_1,\dots, a_k),\ \underline b= (b_1,\dots ,b_\ell)$. This defines the following objects of $\mathrm{Fund}(\g)$:
\[V^{\underline a} := \bigwedge^{a_1}V\otimes \cdots \otimes \bigwedge^{a_k}V,\]
\[V^{\underline b} := \bigwedge^{b_1}V\otimes \cdots \otimes \bigwedge^{b_\ell}V.\]

The goal of this section is to find a diagrammatic basis of webs for $\hom_{\g}\left(V^{\underline a}, V^{\underline b}\right)$. We begin with a dimension count of this space.

\begin{lemma}
\label{lem:dim}
    \[\dim\hom_{\g}(V^{\underline a}, V^{\underline b}) = \begin{cases}
    \binom{k+\ell-2}{k-1} = \binom{k+\ell-2}{\ell-1} &\text{if } |a| = |b|, \\
    0 &\text{otherwise.}
    \end{cases}\]
\end{lemma}

\begin{proof}
    The identity matrix $h\in \gl(1|1)$ acts as the scalar $|\a| = \sum_{i=1}^k a_i$ on $V^{\a}$. Therefore $\hom_{\g}(V^{\a}, V^{\b}) = 0$ unless $|\a| = |\b|$. Hence, suppose $|\a| = |\b|$.
    
    We first claim that 
    \[V^{\underline a} \simeq \bigoplus_{i=0}^{k-1} L(k-i\ \big\vert\  |\a|-k+i)^{\oplus \binom{k-1}{i}}\simeq \bigoplus_{i\ge0} L(i+1\ \big\vert\ |\a|-i-1)^{\oplus \binom{k-1}{i}}.\]

    This follows by an easy induction using the fact that
    \[\bigwedge^mV\otimes \bigwedge^n V = L(1|m-1)\otimes L(1|n-1)\simeq L(1|m+n-1)\oplus L(2|m+n-2).\]
    Writing a similar decomposition of $V^{\underline b}$ into irreducibles and applying Schur's lemma yields

    \begin{align*}
        \dim\hom_\g(V^{\underline a}, V^{\underline b}) &= \sum_{i\ge0} \binom{k-1}{i}\binom{\ell-1}{i}\\
        &= \sum_{i\ge 0}\binom{k-1}{i}\binom{\ell-1}{\ell-i-1}\\
        &= \binom{k+\ell-2}{k-1}. 
    \end{align*}
\end{proof}

\subsection{Diagrammatics}

Fix tuples $\a = (a_1,\dots, a_k), \b=(b_1, \dots, b_\ell)$. A \emph{web} from $\a$ to $\b$ is a morphism in $\infty{-\mathsf{Web}^f}$ from $\a$ to $\b$, drawn using multivalent vertices instead of trivalent ones, such that the graph is \emph{merge-split bipartite} (see \Cref{rem:simplify}). The set of all webs (up to relative height preserving isotopy) from $\a$ to $\b$ is denoted $\W(\a,\b)$.

Just as in the definition of morphisms in $\infty\mathsf{-Web}^f$, we consider webs up to isotopy preserving relative heights of the internal vertices.

Given a web $W$, we get an invariant homomorphism $[W] := \alpha(W)\in \hom_\g(V^{\a}, V^{\b})$ via the realization functor $\alpha$. By the fullness of $\alpha$, the web invariants $\{[W]: W\in \W(\a,\b)\}$ span $\hom_\g(V^{\a},V^{\b})$. We wish to find a natural subset of web invariants which form a basis.

\subsection{Binary flows and monomial expansions} We begin with an alternative way of describing the web invariant $[W]$ in this subsection.

\begin{definition}
    Let $W$ be a web. A \emph{binary flow} (or simply \emph{flow}) $F$ of $W$ is a subgraph of $W$ such that each connected component of $F$ is a directed path. We denote the set of all binary flows of $W$ by $\F_W$.
\end{definition}

\begin{remark}
    We note that this definition is independent of the edge multiplicities in $W$, and depends only on the underlying unlabeled graph.
\end{remark}

Flows allow us to give a monomial expansion of the web invariant $[W]$.

Recall the homogeneous basis $\{x\wedge y^{n-1}, y^n\}$ for $\bigwedge^n V$. Subsets $A\subset [k]$ index a homogeneous basis of monomials $\m_A = v_1\otimes\cdots \otimes v_k$ for $V^{\a}$, via the correspondence

\[
v_i =
\begin{cases}
    x\wedge y^{a_i-1} &\text{if } i\in A,\\
    y^{a_i} &\text{if } i\notin A.
\end{cases}
\]

Tthe space of all $\C$-linear maps $\hom_\C(V^{\a}, V^{\b})$ has a basis of \emph{elementary linear maps} $\{\textbf{E}_{A, B}:A\subset [k], B\subset [\ell]\}$ defined by
\[\textbf{E}_{A,B}(\textbf{m}_S):=
\begin{cases}
    \textbf{m}_B&\text{if } S=A,\\
    0 &\text{otherwise}.
\end{cases}\]

For each flow $F\in \mathcal F_W$, we get subsets $A(F)\subset[k]$ and $B(F)\subset [\ell]$ by looking at the sets of bottom and top boundary vertices in $F$, respectively. This lets us define $\mathbf{E}_F:= \textbf{E}_{A(F), B(F)}\in\hom_\C(V^{\a}, V^{\b})$.

The coefficient $c_F:= \prod_{v\in W}c_{v, F}$ is defined to be the product of local coefficients at each internal vertex $v\in W$, which are given by the rules below:

\begin{equation*}
    v = 
    \left(\vcenter{\hbox{\includegraphics[height=2.5cm]{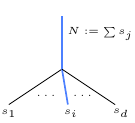}}}\right)
    \implies c_{v,F} = (-1)^{s_1+\cdots + s_{i-1}},
\end{equation*}

\begin{equation*}
    v = 
    \left(\vcenter{\hbox{\includegraphics[height=2.5cm]{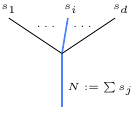}}}\right)
    \implies c_{v,F} = (-1)^{s_1+\cdots + s_{i-1}}\frac{s_i}{N}\binom{N}{s_1,\dots, s_d}.
\end{equation*}

\begin{lemma}\label{prop:flow_expansion}
For each web $W$, the web invariant $[W]$ expands in the basis $\{\textbf{E}_{A, B}\}_{A\subset [k], B\subset [\ell]}$ as
    \[[W]= \sum_{F\in \F_W}c_F\mathbf{E}_F.\]
\end{lemma}

\begin{proof}
    This follows from interpreting $[W]$ as a tensor product and composition of multivalent vertices, and checking the expansion for a single multivalent vertex using the definition of $\alpha$.
\end{proof}

\section{The reduced web basis}\label{sec:rels}

\subsection{Forbidden configurations and reduced webs}
\begin{lemma}\label{lem:red}
    We have the following equalities of web invariants, along with their vertical and horizontal reflections:
    \begin{equation}\label{eq:square_rel1}         \vcenter{\hbox{\includegraphics[width=0.7\linewidth]{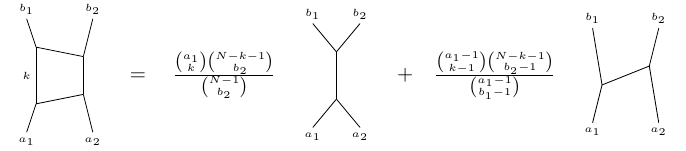}}}
    \quad (\text{for }a_1\ge b_1)
    \end{equation}
    \begin{equation}\label{eq:zipper_rel}
        \vcenter{\hbox{\includegraphics[width=0.9\linewidth]{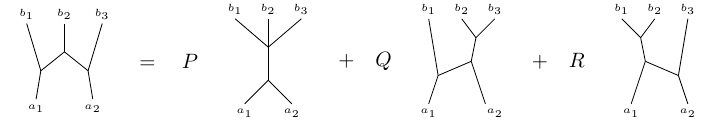}}}
    \end{equation}
    where $N = |\a| = |\b|$, 
    \[P = \binom{N}{b_1, b_2, b_3}^{-1}\Bigg(\binom{a_1}{b_1}\binom{a_2}{b_3} - \frac{b_2+b_3}{b_2}\binom{a_1}{b_1}\binom{a_2-1}{b_3}- \frac{b_1+b_2}{b_2}\binom{a_1-1}{b_1}\binom{a_2}{b_3}\Bigg),\]
    \[Q = \binom{b_2+b_3-1}{b_3}^{-1}\binom{a_2-1}{b_3},\]
    \[R = \binom{b_1+b_2-1}{b_1}^{-1}\binom{a_1-1}{b_1}.\]
\end{lemma}

\begin{proof}
    By \Cref{lem:dim}, the dimension of the morphism space when $k = \ell = 2$ is two, so there must be a linear relation between the three webs in (\ref{eq:square_rel1}). The exact values of the coefficients can be worked out using \Cref{prop:flow_expansion}.

    Similarly, the dimension of the morphism space when $k = 2, \ell = 3$ is three, so there must be a linear relation between the four webs in (\ref{eq:zipper_rel}). The exact values of the coefficients can once again be worked out using \Cref{prop:flow_expansion}.
\end{proof}

\begin{definition}
    We refer to the left hand sides (along with their vertical and horizontal reflections) of (\ref{eq:square_rel1}) and (\ref{eq:zipper_rel}) as \emph{squares} and \emph{zippers}, respectively.
    A web $W$ is \emph{reduced} if it contains no bigons, squares or zippers (in the sense that intersecting $W$ with a rectangular region does not yield a bigon, square, or zipper, even after potentially applying (\ref{eq:assoc}) or (\ref{eq:coassoc})). The set of all reduced webs from $\underline a$ to $\underline b$ is denoted $\RW(\a,\b)$.
\end{definition}

\Cref{lem:red} shows that for any fixed boundary conditions, reduced webs give a spanning set for the corresponding morphism space. The main result of this article is that reduced webs in fact give a basis.

\begin{thm}\label{thm:basis}
    Fix boundary conditions $\underline a$ and $\underline b$ such that $|\underline a| = |\underline b|$. Then $\{[W]: W\in \RW(\a,\b)\}$ is a basis for $\hom_\g(V^{\underline a}, V^{\underline b})$.
\end{thm}

This theorem will follow from \Cref{cor:reduced_web_dim} after we work out some combinatorial properties of reduced webs in the next subsection. As a corollary of \Cref{thm:basis}, we obtain the following alternative presentation of $1|1\mathsf{-Web}$.

\begin{corollary}\label{cor:alt_presentation}
    Let $1|1\mathsf{-Web}'$ be the quotient of $\infty\mathsf{-Web}^f$ by the relations (\ref{eq:assoc}), (\ref{eq:coassoc}), (\ref{eq:bigon}), (\ref{eq:square_rel1}), and (\ref{eq:zipper_rel}), along with the vertical and horizontal reflections of the last two relations. Then $\alpha$ descends to an equivalence of categories between $1|1\mathsf{-Web}'$ and $\mathrm{Fund}(\g)$.
\end{corollary}

\begin{remark}
    While it is not difficult to derive (\ref{eq:funny}) from (\ref{eq:square_rel1}), it is not clear to us how to derive (\ref{eq:square_rel1}) and (\ref{eq:zipper_rel}) from the presentation in \Cref{thm:tvw}. 
\end{remark}

In fact, the reduced web basis comes equipped with an algorithm to decompose \emph{any} web into a linear combination of reduced webs.

\begin{algorithm}[Reduction algorithm]\label{algo:red} Let $W$ be any web.
\begin{enumerate}
    \item Resolve all bigons using (\ref{eq:bigon}), and make $W$ merge-split bipartite using multivalent vertices via (\ref{eq:assoc}) and (\ref{eq:coassoc}) (c.f. \Cref{rem:simplify}).
    \item Fix a face of $W$. 
    \begin{itemize}
        \item If the chosen face has $\ge6$ edges, then resolve the zipper at the top vertex using (\ref{eq:zipper_rel}), and then introduce more multivalent vertices to make $W$ merge-split bipartite. This reduces the number of edges in the chosen face.
        \item If the chosen face is a square, resolve it using (\ref{eq:square_rel1}).
    \end{itemize}
    Repeat for all other faces of $W$ (c.f. proof of \Cref{prop:acyclic}).
    \item Resolve all zippers using (\ref{eq:zipper_rel}), after possibly replacing some multivalent vertices with trivalent ones.
    \item Reintroduce multivalent vertices to make all resultant webs merge-split bipartite.
\end{enumerate}
\end{algorithm}

To see that this algorithm terminates, observe that:
\begin{itemize}
    \item each application of step $(2)$ strictly decreases the number of faces of $W$, and
    \item each application of (\ref{eq:zipper_rel}) or its vertical reflection strictly decreases the number of directed edges from a split to a merge vertex, while not reintroducing any faces.
\end{itemize}

\subsection{Combinatorics of Reduced webs}\label{sec:proof}

We begin by showing that reduced webs contain no cycles (as undirected graphs), i.e., they are forests.

\begin{figure}
    \centering
    \includegraphics[height=6cm]{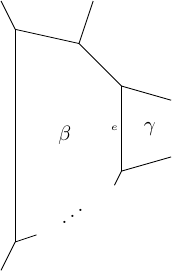}
    \caption{The situation described in the last paragraph of the proof of \Cref{prop:acyclic}.}
    \label{fig:acyclicity_proof}
\end{figure}

\begin{proposition}\label{prop:acyclic}
    Let $W$ be a reduced web. Then the underlying undirected graph of $W$ contains no cycles.
\end{proposition}
\begin{proof}
    Suppose for the sake of contradiction that $\beta$ is a face of $W$. Then the boundary of $\beta$ consists of two directed paths from the bottom vertex to the top vertex. Let the lengths of the left and right boundary paths be $l_\beta$ and $r_\beta$ respectively. If $l_\beta>1$ and $r_\beta>1$, then the top vertex of $\beta$ is a zipper, contradicting the assumption that $W$ is reduced. Thus $l_\beta = 1$ or $r_\beta=1$.

    By horizontally flipping the web if necessary, we may assume that $l_\beta = 1$. Since $W$ is a reduced web, $\beta$ cannot be a square, so $r_\beta\ge 4$. Let $e$ be third edge from the top in the right path. Then the top vertex of $e$ is a merge vertex, different from the top of $\beta$. Similarly, the bottom of $e$ is a split vertex, different from the bottom of $\beta$. By planarity, this forces $W$ to contain a face $\gamma$ whose left boundary is the single edge $e$. In summary, given $\beta$ with $l_\beta = 1$, we have shown that there is a face $\gamma$ to the right of $\beta$ with $l_{\gamma} = 1$ (see \Cref{fig:acyclicity_proof} for a depiction of this situation). Repeating the argument gives an infinite chain of faces, contradicting the finiteness of $W$.
\end{proof}

Next, we show that for boundary conditions $\underline a = (a_1,\dots, a_k)$ and $\underline b = (\b_1,\dots, b_\ell)$, the number of reduced webs is $|\RW(\a,\b)| =\binom{k+\ell-2}{k-1}$. Our strategy will be to produce an explicit bijection to a set of the desired cardinality.

Let
\[S(k, \ell) = \{(A,B)\in 2^{[k-1]}\times 2^{[\ell-1]} :  |A| = |B| \}.\]
In words, $S(k, \ell)$ is the set of pairs of subsets of $[k-1]$ and $[\ell-1]$ of the same cardinality. Clearly,

\[|S(k,\ell)| = \sum_{i\ge0}\binom{k-1}{i}\binom{\ell-1}{i} = \sum_{i\ge0}\binom{k-1}{k-i-1}\binom{\ell-1}{i} = \binom{k+\ell-2}{k-1}.\]

We now construct a map $\Phi_{\a, \b}: \RW(\a,\b)\to S(k, \ell)$. Given $W\in \RW(\a,\b)$, consider the \emph{rightmost maximal flow} $F_W\in \F_W$, and let $\Phi_{\a,\b}(W)$ be given by the endpoints of $F_W$.

More precisely, $F_W$ is built up inductively as follows. Declare the last component $f$ of $F_W$ to be the flow that starts at $a_k$ and goes maximally right through $W$ to $b_\ell$. Let $W'$ be the induced subgraph of $W$ consisting of all vertices not in $f$. Then $W'$ is also reduced, and we recursively define $F_W := F_{W'}\sqcup f$. Once $F_W$ is defined,
\[A_W = \{i\in [k-1]: a_i \text{ is a source for }F_W\},\]
\[B_W = \{j\in [\ell-1]: b_j \text{ is a sink for }F_W\},\]
defines the element $\Phi_{\a, \b}(W):= (A_W,B_W)\in S(k,\ell)$.\\

\begin{figure}[tph]
    \centering
    \includegraphics[height=5cm]{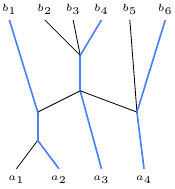}
    \caption{The rightmost maximal flow $F_W$ on the web $W$, shown in blue. The endpoints of this flow show that $\Phi_{\a,\b}(W) = (\{2,3,4\}, \{1,4,6\})$.}
    \label{fig:rmflow}
\end{figure}

\begin{figure}
    \centering
    \includegraphics[height=4.5cm]{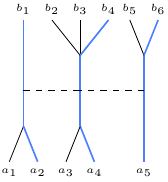}
    \caption{The image under $\Psi_{\a,\b}$ of $(\{2, 4\}, \{1, 4\})$. The orientation of the left horizontal dashed edge depends on whether $a_1+a_2<b_1$ or $a_1+a_2>b_1$ (and this edge vanishes if $a_1+a_2=b_1$). Similarly, the orientation of the right dashed edge depends on $a_5$ versus $b_5+b_6$ (or equivalently, $a_1+a_2+a_3+a_4$ versus $b_1+b_2+b_3+b_4$). Lastly, the relative heights of these edges (in case they are both oriented in the same direction) are chosen in the unique way so that there is no zipper.}
    \label{fig:psi}
\end{figure}

In the other direction, we construct $\Psi_{\a, \b}: S(k,\ell)\to \RW(\a,\b)$ as follows. Let $(A,B) \in S(k,\ell)$ with $A = \{i_1<\cdots< i_t\}$ and $B= \{j_1<\cdots < j_t\}$. For each $1\le s\le t$, draw the graph which is a merge of vertices $i_s,\dots, i_{s+1}-1$ in $\a$, followed by a split to vertices $j_s,\dots, j_{s+1}-1$ of $\b$ (with the convention that $i_{t+1}-1 = k$ and $j_{t+1}-1 = \ell$). As defined, this is not a web, since the algebraic sum of edge multiplicities at an internal vertex is not $0$. However, adding a single edge between each of these components in order to satisfy the algebraic sum condition yields a genuine web from $\a$ to $\b$. This is a reduced web, which we define to be $\Psi_{\a,\b}(A,B)$. See \Cref{fig:psi} for an example of this construction.

It is clear from the construction that $\Phi_{\a,\b}\o \Psi_{\a,\b} = \mathrm{id}_{S(k,\ell)}$.

\begin{remark}
    $\Phi_{\a,\b}\o\Psi_{\a,\b} = \mathrm{id}_{S(k,\ell)}$ implies that $\Psi_{\a,\b}$ is injective, so that $|\mathrm{im}(\Psi_{\a,\b})| = |S(k,\ell)| = |\dim\hom_{\g}(V^{\a},V^{\b})|$. Moreover, one can use a leading term argument with respect to an appropriate term order on the monomial basis to show that $\{[W]: W\in \mathrm{im}(\Psi_{\a,\b})\}$ is linearly independent. This suffices to show that $\mathrm{im}(\Psi_{\a,\b})$ gives a web basis for $\hom_\g(V^{\a}, V^{\b})$. However, in order to obtain the reduction algorithm, we will show that $\Psi_{\a,\b}$ is a bijection so that $\mathrm{im}(\Psi_{\a,\b}) = \RW(\a,\b)$.
\end{remark}

\begin{proposition}
    $\Phi_{\a, \b}$ is injective.
\end{proposition}

\begin{proof}
Suppose $W\in \RW(\a,\b)$, with $\Phi_{\a,\b}(W) = (A_W, B_W)\in S(k,\ell)$. We wish to show that $W = \Psi_{\a,\b}(A_W, B_W)$. We proceed by induction on $|A_W| = |B_W|$. 
    
    Recall that $F_W\in \mathcal{F}_W$ is the rightmost maximal binary flow in $W$. We claim that the rightmost component $f$ of $F_W$ cannot contain a split followed by a merge. If it did, planarity would force a face in $W$, contradicting \Cref{prop:acyclic}. In particular, $f$ contains at most two internal vertices, with a merge followed by a split when equality is achieved. 
    
    If $i_0 = \max A_W$, then the boundary vertices $a_{i_0+1},\dots, a_k$ must all have edges to the merge vertex in $f$. Similarly, if $j_0=\max B_W$, then the boundary vertices $b_{j_0+1},\dots, b_\ell$ must all have edges to the split vertex in $f$. Moreover, no other boundary vertex has an edge connecting to $f$, so any other edges of $W$ that touch $f$ must connect to the second rightmost component $f'$ of $F_W$. Once again, by the acyclicity from \Cref{prop:acyclic}, there is at most one edge $e$ of $W$ that connects $f$ and $f'$.

    The orientation and multiplicity of $e$ are completely determined by the difference
    \[\sum_{i=i_0}^k a_i - \sum_{j=j_0}^\ell b_j.\]

    In summary, the edges that touch $f$ are those connecting to the boundary vertices $a_{i_0},\dots, a_k, b_{j_0},\dots, b_\ell$, and potentially the edge $e$ connecting $f$ and $f'$. Deleting $f$ and all these edges gives a smaller reduced web $W'$.
    
    Let the boundary conditions of $W'$ be $\a',\b'$, with $\Phi_{\a',\b'}(W') = (A_{W'}, B_{W'}) = \big(A_W\setminus \{\max A_W\}, B_W\setminus \{\max B_W\}\big)$. By the inductive hypothesis, $W' = \Psi_{\a',\b'}(A_{W'}, B_{W'})$. Therefore, by the definition of $\Psi_{\a, \b}$ it follows that $W = \Psi_{\a,\b}(A_W, B_W)$, as desired.  
\end{proof}

\begin{corollary}\label{cor:reduced_web_dim}
    $\Phi_{\a,\b}$ and $\Psi_{\a, \b}$ give inverse bijections between $\RW(\a,\b)$ and $S(k, \ell)$. In particular,
    \[|\RW(\a,\b)| = |S(k, \ell)| = \binom{k+\ell-2}{k-1}.\]
\end{corollary}

\begin{proof}[Proof of \Cref{thm:basis}]
    $\{[W]: W\in\RW(\a,\b)\}$ spans $\hom_{\g}(V^{\a}, V^{\b})$ by \Cref{lem:red}, and \Cref{cor:reduced_web_dim} shows that its cardinality is equal to $\dim\hom_\g(V^{\a}, V^{\b})$ by \Cref{lem:dim}. This proves the theorem.
\end{proof}

\section{Clasped webs}\label{sec:clasped-webs}

\subsection{Clasped webs}

Inside $\bigwedge^{a_1}V\otimes \cdots \otimes \bigwedge^{a_i}V$, the tensor product of highest weight vectors generates an irreducible submodule, the \emph{Cartan component}, isomorphic to $$L(i\ \big\vert\ |\a| - i).$$ We depict this representation diagrammatically by ``clasping" together $i$ boundary vertices of multiplicities $a_1, \dots, a_i$ from left to right.

\begin{definition}
    A \emph{clasped boundary condition} is a pair of tuples $\underline c = (c_1, \dots, c_k)$ and $\underline a = (a_1, \dots, a_{|\underline c|})$, denoted by $(\a;\underline c)$. This defines a tensor product of irreducible representations 
\[L(\underline a; \underline c)\subset V^{\a}\]
by clasping together the first $c_1$ vertices, the next $c_2$ vertices, and so on in $V^{\a}$.
\end{definition}

\begin{figure}[tph]
    \centering
    \includegraphics[width=0.7\linewidth]{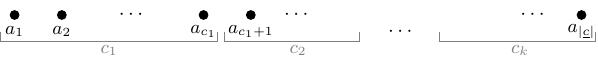}
    \caption{The clasped boundary condition $(\a;\underline c)$ with $k$ clasps (drawn in gray).}
    \label{fig:clasp}
\end{figure}

Pictorially, $k$ is the number of clasps and $c_i$ is the number of boundary vertices in the $i^\text{th}$ clasp, so that $|\underline c|$ is the total number of boundary vertices. The tuple $\a$ just records the labels of the boundary vertices, as before. 

The tensor product of irreducible representations defined above is isomorphic to the tensor product
\[L(\underline a; \underline c)\simeq L(c_1|p_1)\otimes \cdots \otimes L(c_k|p_k),\]
where $p_i = \big(\sum_{j=c_1+\cdots +c_{i-1}+1}^{c_1+\cdots+c_i}a_j\big) - c_i$, so that $|\underline c|+|\underline p| = |\a|$.

Similarly, let $\underline d = (d_1, \dots, d_\ell)$ and $\b = (b_1, \dots, b_{|\underline d|})$ be another clasped boundary condition, so that 
\[L(\b; \underline d)\subset V^{\b}.\]

We then have a projection map
\[\proj:\hom_\g(V^{\a}, V^{\b})\to \hom_\g\big(L(\a;\underline c), L(\b;\underline d)\big),\]
given by the restriction to $L(\a;\underline c)\subset V^{\a}$ followed by the projection to the direct summand $L(\b,\underline d)\subset V^{\b}$.

The goal of this section is to describe $\proj$ in terms of the web basis for $\hom_\g(V^{\a}, V^{\b})$. We depict the image of a web invariant $[W]\in \hom_\g(V^{\a}, V^{\b})$ under $\pi_{(\a;\underline c), (\b;\underline d)}$ diagrammatically by introducing ``external clasps" around the boundary vertices of $W$, and call the resultant a \emph{clasped web}.

\begin{figure}[tph]
    \centering
    \includegraphics[height=4cm]{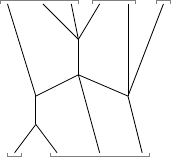}
    \caption{Example of a clasped web with $\underline c = (1, 3)$ and $\underline d = (3, 2, 1)$, where the external clasps are drawn in gray.}
    \label{fig:clasped_web}
\end{figure}

\subsection{Dimensions of clasped homomorphism spaces}
First, we record the dimension of the clasped homomorphism space, generalizing \Cref{lem:dim}.

\begin{lemma}\label{lem:dim_clasped}
    \[\dim \hom_\g\big(L(\a; \underline c), L(\b;\underline d)\big) = 
    \begin{cases}
        \binom{k+\ell-2}{k-1+|\underline d|-|\underline c|} = \binom{k+\ell-2}{\ell-1 +|\underline c|-|\underline d|} & \text{ if } |\a|= |\b|,\\
        0 &\text{ otherwise.}
    \end{cases}
    \]
\end{lemma}

\begin{proof}
    Writing 
    \(L(\a; \underline c) \simeq L(c_1|p_1)\otimes\cdots \otimes L(c_k|p_k)\)
    as before, we claim that 
    \[L(\a;\underline c)\simeq \bigoplus_{i=0}^{k-1}L(|\underline c|-i\ \big\vert\ |\underline p|+i)^{\oplus \binom{k-1}{i}}.\]
    This follows from an easy induction using the fact that for any $(u_1|v_1), (u_2|v_2)\in \Lambda^{\mathrm{pol}}\setminus \{(0|0)\}$,
    \[L(u_1|v_1)\otimes L(u_2|v_2) \simeq L(u_1+u_2|v_1+v_2)\oplus L(u_1+u_2-1| v_1+v_2+1).\]
    Similarly, writing $L(\b;\underline d) \simeq L(b_1|q_1)\otimes \cdots \otimes L(b_\ell|q_\ell)$, we find that
    \[L(\b; \underline d) \simeq \bigoplus_{j=0}^{\ell-1} L(|\underline d|-j\ \big\vert\ |\underline q|+j)^{\oplus \binom{\ell-1}{j}}. \]
    By Schur's lemma, we see that if $|\a|\ne |\b|$, then $\dim\hom_{\g}\big(L(\a;\underline c), L(\b;\underline d)\big) = 0$, and if $|\a| = |\b|$, then
    \begin{align*}
        \dim\hom_\g\big(L(\a;\underline c), L(\b;\underline d)\big) &= \sum_{i,j\ge 0}\delta_{|\underline c|-i, |\underline d|-j}\binom{k-1}{i}\binom{\ell-1}{j}\\
        &= \sum_{i\ge 0}\binom{k-1}{i}\binom{\ell-1}{i -|\underline c|+|\underline d|}\\
        &= \sum_{i\ge0}\binom{k-1}{i}\binom{\ell-1}{\ell-1+|\underline c|-|\underline d|-i}\\
        &= \binom{k+\ell-2}{\ell-1+|\underline c|-|\underline d|},
    \end{align*}
    where $\delta_{|\underline c|-i, |\underline d|-j}$ in the first line is the Kronecker delta.
\end{proof}

\subsection{Non-convex webs}

\begin{definition}
    Fix clasped boundary conditions $(\a;\underline c)$ and $(\b;\underline d)$. A reduced web $W\in \RW(\a,\b)$ is \emph{non-convex} if no two boundary vertices within the same clasp share a common neighbor in $W$. We say that $W\in \RW(\a,\b)$ is \emph{partially convex} if it is not non-convex.
\end{definition}

\begin{figure}[tph]
    \centering
    \includegraphics[height=2.5cm]{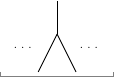}
    \caption{A web is \emph{partially convex} if and only if it contains the above local configuration (or its vertical reflection) within a clasp.}
    \label{fig:partially_convex}
\end{figure}

For example, the web in \Cref{fig:clasped_web} is partially convex because the leftmost clasp at the top boundary violates the non-convexity condition.
The following lemma motivates the definition of non-convexity. 

\begin{lemma}\label{lem:non-convex}
    If $W\in \RW(\a,\b)$ is partially convex, then $\proj([W]) = 0$.
\end{lemma}

\begin{proof}
    This can be checked directly using \Cref{prop:flow_expansion}, noting that there are no flows incident to both boundary vertices in a partially convex local configuration as in \Cref{fig:partially_convex}.
\end{proof}

\begin{corollary}\label{cor:spanning}
    A spanning set for $\hom_\g\big(L(\a;\underline c), (\b;\underline d)\big)$ is given by
    \[\{\proj([W]): W\in \RW(\a,\b) \text{ non-convex}\}.\]
\end{corollary}

The main result of this section is that the projections of the non-convex webs in fact form a basis for $\hom_\g\big(L(\a;\underline c), (\b;\underline d)\big)$. This follows from a dimension counting argument.

\begin{lemma}\label{lem:clasped_web_count}
    Fix clasped boundary conditions $(\a;\underline c)$ and $(\b;\underline d)$ such that $|\a| = |\b|$. Then
    \[|\{W\in \RW(\a,\b): W \text{ is non-convex}\}| = \binom{k+\ell-2}{k-1+|\underline d|-|\underline c|} = \binom{k+\ell-2}{\ell-1 +|\underline c|-|\underline d|}.\]
\end{lemma}

\begin{proof}
    We will exploit the explicit description of $\RW(\a,\b)$ as the image of $S(k,\ell)$ under $\Psi_{\a,\b}$. In particular, if $W\in \RW(\a,\b)$ has rightmost maximal flow $F_W$ with endpoints $(A_W, B_W)\in S(k,\ell)$, then $W = \Psi_{\a,\b}(A_W,B_W)$. It is clear from this description that $W$ is non-convex if and only if for each clasp, if its length is $t$, then the leftmost $t-1$ vertices of the clasp are contained in $F_W$. Thus, the only data left to record is which clasps also have their rightmost vertex contained in $F_W$. We say a clasp is \emph{marked} if its rightmost vertex is contained in $F_W$.
    
    Obviously, the rightmost clasps at the bottom and the top boundary are marked. Say we pick $i$ out of the remaining bottom $k-1$ clasps and $j$ out of the remaining top $\ell-1$ clasps to be marked. Counting the number of components of $F_W$, we see that
    \[|\underline c|-k +i+1 = |\underline d|-\ell+j+1,\]
    and hence $j = i-k+\ell+|\underline c|-|\underline d|$. For ease of notation, let $\Delta = |\underline c|-|\underline d|$, so that $j=i-k+\ell + \Delta$. Since we may choose the marked clasps freely, we see that the number of non-convex webs is
    \begin{align*}
        \sum_{i\ge0}\binom{k-1}{i}\binom{\ell-1}{i-k+\ell+\Delta} &= \sum_{i\ge0}\binom{k-1}{i}\binom{\ell-1}{k-\Delta-1-i}\\
        &= \binom{k+\ell-2}{k-\Delta-1}\\
        &= \binom{k+\ell-2}{k-1+|\underline d|-|\underline c|}.
    \end{align*}    
\end{proof}

This gives us the main theorem of this section.

\begin{thm}\label{thm:clasped}
    Fix clasped boundary conditions $(\a;\underline c)$ and $(\b;\underline d)$. Then $\proj$ sends the partially convex webs in the web basis $\{[W]: W\in \RW(\a,\b)\}$ to $0$, and the non-convex webs to a basis of $\hom_\g\big(L(\a;\underline c), (\b;\underline d)\big)$.
\end{thm}

\begin{proof}
    The first assertion is \Cref{lem:non-convex}. The second assertion follows from \Cref{cor:spanning} after comparing dimensions using \Cref{lem:dim_clasped} and \Cref{lem:clasped_web_count}.
\end{proof}

\section*{Acknowledgements}
I wish to thank my advisors Christian Gaetz and Vera Serganova for many helpful conversations. 

\section*{AI disclosure statement}
There was no use of AI in the preparation of this paper, and the author is completely responsible for its results and proofs.

\bibliographystyle{amsalphavar}
\bibliography{main}

@article {Douglas-Kenyon-Shi,
    AUTHOR = {Douglas, Daniel C. and Kenyon, Richard and Shi, Haolin},
     TITLE = {Dimers, webs, and local systems},
   JOURNAL = {Trans. Amer. Math. Soc.},
  FJOURNAL = {Transactions of the American Mathematical Society},
    VOLUME = {377},
      YEAR = {2024},
    NUMBER = {2},
     PAGES = {921--950},
      ISSN = {0002-9947,1088-6850},
   MRCLASS = {82B20 (57M15)},
  MRNUMBER = {4688540},
       DOI = {10.1090/tran/8993},
       URL = {https://doi.org/10.1090/tran/8993},
}

@article {Cautis-Kamnitzer-Morrison,
    AUTHOR = {Cautis, Sabin and Kamnitzer, Joel and Morrison, Scott},
     TITLE = {Webs and quantum skew {H}owe duality},
   JOURNAL = {Math. Ann.},
  FJOURNAL = {Mathematische Annalen},
    VOLUME = {360},
      YEAR = {2014},
    NUMBER = {1-2},
     PAGES = {351--390},
      ISSN = {0025-5831,1432-1807},
   MRCLASS = {17B37 (17B10 17B20)},
  MRNUMBER = {3263166},
MRREVIEWER = {Iwan\ Praton},
       DOI = {10.1007/s00208-013-0984-4},
       URL = {https://doi.org/10.1007/s00208-013-0984-4},
}

@article {Fomin-Pylyavskyy-advances,
    AUTHOR = {Fomin, Sergey and Pylyavskyy, Pavlo},
     TITLE = {Tensor diagrams and cluster algebras},
   JOURNAL = {Adv. Math.},
  FJOURNAL = {Advances in Mathematics},
    VOLUME = {300},
      YEAR = {2016},
     PAGES = {717--787},
      ISSN = {0001-8708,1090-2082},
   MRCLASS = {13F60 (05E99 13A50 15A72 20G05)},
  MRNUMBER = {3534844},
MRREVIEWER = {Olga\ Kravchenko},
       DOI = {10.1016/j.aim.2016.03.030},
       URL = {https://doi.org/10.1016/j.aim.2016.03.030},
}

@article {Fontaine,
    AUTHOR = {Fontaine, Bruce},
     TITLE = {Generating basis webs for {${\rm SL_n}$}},
   JOURNAL = {Adv. Math.},
  FJOURNAL = {Advances in Mathematics},
    VOLUME = {229},
      YEAR = {2012},
    NUMBER = {5},
     PAGES = {2792--2817},
      ISSN = {0001-8708,1090-2082},
   MRCLASS = {20G15 (20G05)},
  MRNUMBER = {2889146},
MRREVIEWER = {Timothy\ C.\ Burness},
       DOI = {10.1016/j.aim.2012.01.016},
       URL = {https://doi.org/10.1016/j.aim.2012.01.016},
}

@article {Khovanov,
    AUTHOR = {Khovanov, Mikhail},
     TITLE = {sl(3) link homology},
   JOURNAL = {Algebr. Geom. Topol.},
  FJOURNAL = {Algebraic \& Geometric Topology},
    VOLUME = {4},
      YEAR = {2004},
     PAGES = {1045--1081},
      ISSN = {1472-2747,1472-2739},
   MRCLASS = {57M27 (18G60 57R56)},
  MRNUMBER = {2100691},
MRREVIEWER = {Justin\ Sawon},
       DOI = {10.2140/agt.2004.4.1045},
       URL = {https://doi.org/10.2140/agt.2004.4.1045},
}

@article {Kuperberg,
    AUTHOR = {Kuperberg, Greg},
     TITLE = {Spiders for rank {$2$} {L}ie algebras},
   JOURNAL = {Comm. Math. Phys.},
  FJOURNAL = {Communications in Mathematical Physics},
    VOLUME = {180},
      YEAR = {1996},
    NUMBER = {1},
     PAGES = {109--151},
      ISSN = {0010-3616,1432-0916},
   MRCLASS = {17B10 (22E60 81R05)},
  MRNUMBER = {1403861},
MRREVIEWER = {Stefano\ Capparelli},
       URL = {http://projecteuclid.org/euclid.cmp/1104287237},
}

@article {Petersen-Pylyavskyy-Rhoades,
    AUTHOR = {Petersen, T. Kyle and Pylyavskyy, Pavlo and Rhoades, Brendon},
     TITLE = {Promotion and cyclic sieving via webs},
   JOURNAL = {J. Algebraic Combin.},
  FJOURNAL = {Journal of Algebraic Combinatorics. An International Journal},
    VOLUME = {30},
      YEAR = {2009},
    NUMBER = {1},
     PAGES = {19--41},
      ISSN = {0925-9899,1572-9192},
   MRCLASS = {05E10},
  MRNUMBER = {2519848},
MRREVIEWER = {Gregory\ S.\ Warrington},
       DOI = {10.1007/s10801-008-0150-3},
       URL = {https://doi.org/10.1007/s10801-008-0150-3},
}

@article {Reshetikhin-Turaev,
    AUTHOR = {Reshetikhin, N. Yu. and Turaev, V. G.},
     TITLE = {Ribbon graphs and their invariants derived from quantum
              groups},
   JOURNAL = {Comm. Math. Phys.},
  FJOURNAL = {Communications in Mathematical Physics},
    VOLUME = {127},
      YEAR = {1990},
    NUMBER = {1},
     PAGES = {1--26},
      ISSN = {0010-3616,1432-0916},
   MRCLASS = {57M25 (16W30 17B35)},
  MRNUMBER = {1036112},
MRREVIEWER = {Louis\ H.\ Kauffman},
       URL = {http://projecteuclid.org/euclid.cmp/1104180037},
}

@article {GPPSS-sl4,
    AUTHOR = {Gaetz, Christian and Pechenik, Oliver and Pfannerer, Stephan
              and Striker, Jessica and Swanson, Joshua P.},
     TITLE = {Rotation-invariant web bases from hourglass plabic graphs},
   JOURNAL = {Invent. Math.},
  FJOURNAL = {Inventiones Mathematicae},
    VOLUME = {243},
      YEAR = {2026},
    NUMBER = {2},
     PAGES = {703--804},
      ISSN = {0020-9910,1432-1297},
   MRCLASS = {05E10 (13F60 17B37 18M30)},
  MRNUMBER = {5008156},
MRREVIEWER = {Teng\ Fang},
       DOI = {10.1007/s00222-025-01385-2},
       URL = {https://doi.org/10.1007/s00222-025-01385-2},
}

@article{Elias,
  doi = {10.48550/ARXIV.1510.06840},
  url = {https://arxiv.org/abs/1510.06840},
  author = {Elias, Ben},
  title = {Light ladders and clasp conjectures},
  publisher = {arXiv},
  journal={preprint},
  year = {2015},
  note={\url{arxiv:1510.06840}},
  pages={69 pages}
}

@article {Le.Sikora,
    AUTHOR = {L\^e, Thang T. Q. and Sikora, Adam S.},
     TITLE = {Stated {${\rm SL}(n)$}-skein modules and algebras},
   JOURNAL = {J. Topol.},
  FJOURNAL = {Journal of Topology},
    VOLUME = {17},
      YEAR = {2024},
    NUMBER = {3},
     PAGES = {Paper No. e12350, 93},
      ISSN = {1753-8416,1753-8424},
   MRCLASS = {57K31 (17B37 20G42 57K16)},
  MRNUMBER = {4822918},
MRREVIEWER = {Uwe\ Kaiser},
       DOI = {10.1112/topo.12350},
       URL = {https://doi.org/10.1112/topo.12350},
}

@article {grant16,
    AUTHOR = {Grant, Jonathan},
     TITLE = {A generators and relations description of a representation
              category of {$U_q(\mathfrak{gl}(1|1))$}},
   JOURNAL = {Algebr. Geom. Topol.},
  FJOURNAL = {Algebraic \& Geometric Topology},
    VOLUME = {16},
      YEAR = {2016},
    NUMBER = {1},
     PAGES = {509--539},
      ISSN = {1472-2747,1472-2739},
   MRCLASS = {17B37 (57M25)},
  MRNUMBER = {3470707},
MRREVIEWER = {Pedro\ Vaz},
       DOI = {10.2140/agt.2016.16.509},
       URL = {https://doi.org/10.2140/agt.2016.16.509},
}

@article {Sartori2013CategorificationOT,
    AUTHOR = {Sartori, Antonio},
     TITLE = {Categorification of tensor powers of the vector representation
              of {$U_q(\mathfrak{gl}(1|1))$}},
   JOURNAL = {Selecta Math. (N.S.)},
  FJOURNAL = {Selecta Mathematica. New Series},
    VOLUME = {22},
      YEAR = {2016},
    NUMBER = {2},
     PAGES = {669--734},
      ISSN = {1022-1824,1420-9020},
   MRCLASS = {17B10 (17B37 20C08)},
  MRNUMBER = {3477333},
MRREVIEWER = {Simon\ Riche},
       DOI = {10.1007/s00029-015-0202-1},
       URL = {https://doi.org/10.1007/s00029-015-0202-1},
}

@article {tubbenhauer-vaz-wedrich,
    AUTHOR = {Tubbenhauer, Daniel and Vaz, Pedro and Wedrich, Paul},
     TITLE = {Super {$q$}-{H}owe duality and web categories},
   JOURNAL = {Algebr. Geom. Topol.},
  FJOURNAL = {Algebraic \& Geometric Topology},
    VOLUME = {17},
      YEAR = {2017},
    NUMBER = {6},
     PAGES = {3703--3749},
      ISSN = {1472-2747,1472-2739},
   MRCLASS = {57M25 (81R50)},
  MRNUMBER = {3709658},
MRREVIEWER = {Benjamin\ Cooper},
       DOI = {10.2140/agt.2017.17.3703},
       URL = {https://doi.org/10.2140/agt.2017.17.3703},
}

@article{enugandla-gaetz,
  doi = {10.48550/ARXIV.2512.08817},
  url = {https://arxiv.org/abs/2512.08817},
  author = {Pranav Enugandla and Christian Gaetz},
  title = {Clasped web bases from hourglass plabic graphs},
  publisher = {arXiv},
  journal={preprint},
  year = {2025},
  note={\url{arxiv:2512.08817}},
  pages={40 pages}
}

@article {queffelec-sartori,
    AUTHOR = {Queffelec, Hoel and Sartori, Antonio},
     TITLE = {Mixed quantum skew {H}owe duality and link invariants of type
              {$A$}},
   JOURNAL = {J. Pure Appl. Algebra},
  FJOURNAL = {Journal of Pure and Applied Algebra},
    VOLUME = {223},
      YEAR = {2019},
    NUMBER = {7},
     PAGES = {2733--2779},
      ISSN = {0022-4049,1873-1376},
   MRCLASS = {17B37 (17B35 20G43 57M25 57M27 81R50)},
  MRNUMBER = {3912946},
MRREVIEWER = {Li-Xia\ Ye},
       DOI = {10.1016/j.jpaa.2018.09.014},
       URL = {https://doi.org/10.1016/j.jpaa.2018.09.014},
}

@incollection {queffelec2015noteglmnlinkinvariants,
    AUTHOR = {Queffelec, Hoel and Sartori, Antonio},
     TITLE = {A note on {$\mathfrak{gl}_{m|n}$} link invariants and the
              {HOMFLY}-{PT} polynomial},
 BOOKTITLE = {Knots, low-dimensional topology and applications},
    SERIES = {Springer Proc. Math. Stat.},
    VOLUME = {284},
     PAGES = {277--286},
 PUBLISHER = {Springer, Cham},
      YEAR = {2019},
      ISBN = {978-3-030-16031-9; 978-3-030-16030-2},
   MRCLASS = {57M27 (57M25 81R50)},
  MRNUMBER = {3986050},
MRREVIEWER = {Haimiao\ Chen},
       DOI = {10.1007/978-3-030-16031-9\_13},
       URL = {https://doi.org/10.1007/978-3-030-16031-9_13},
}

@article {Sartori2013TheAP,
    AUTHOR = {Sartori, Antonio},
     TITLE = {The {A}lexander polynomial as quantum invariant of links},
   JOURNAL = {Ark. Mat.},
  FJOURNAL = {Arkiv f\"or Matematik},
    VOLUME = {53},
      YEAR = {2015},
    NUMBER = {1},
     PAGES = {177--202},
      ISSN = {0004-2080,1871-2487},
   MRCLASS = {17B37 (57M27)},
  MRNUMBER = {3319619},
MRREVIEWER = {Rinat\ M.\ Kashaev},
       DOI = {10.1007/s11512-014-0196-5},
       URL = {https://doi.org/10.1007/s11512-014-0196-5},
}

@book {cheng-wang,
    AUTHOR = {Cheng, Shun-Jen and Wang, Weiqiang},
     TITLE = {Dualities and representations of {L}ie superalgebras},
    SERIES = {Graduate Studies in Mathematics},
    VOLUME = {144},
 PUBLISHER = {American Mathematical Society, Providence, RI},
      YEAR = {2012},
     PAGES = {xviii+302},
      ISBN = {978-0-8218-9118-6},
   MRCLASS = {17B10},
  MRNUMBER = {3012224},
MRREVIEWER = {Aleksandr\ Nikolaevich\ Sergeev},
       DOI = {10.1090/gsm/144},
       URL = {https://doi.org/10.1090/gsm/144},
}

@book {musson,
    AUTHOR = {Musson, Ian M.},
     TITLE = {Lie superalgebras and enveloping algebras},
    SERIES = {Graduate Studies in Mathematics},
    VOLUME = {131},
 PUBLISHER = {American Mathematical Society, Providence, RI},
      YEAR = {2012},
     PAGES = {xx+488},
      ISBN = {978-0-8218-6867-6},
   MRCLASS = {17-02 (16S30 17B35)},
  MRNUMBER = {2906817},
MRREVIEWER = {Aleksandr\ Nikolaevich\ Sergeev},
       DOI = {10.1090/gsm/131},
       URL = {https://doi.org/10.1090/gsm/131},
}
\end{document}